\documentclass[12pt]{amsart}

\usepackage{times}

\usepackage{graphicx, xspace}

\theoremstyle{plain}
\newtheorem{lemma}{Lemma}
\newtheorem{theorem}[lemma]{Theorem}
\newtheorem{proposition}[lemma]{Proposition}
\newtheorem{problem}[lemma]{Problem}
\newtheorem{conjecture}[lemma]{Conjecture}
\newtheorem{corollary}[lemma]{Corollary}

\theoremstyle{remark}

\newcommand{\normchar}{\tilde{\chi}}
\newcommand{\HH}{\tilde{H}}
\newcommand{\falling}[2]{(#1)_{#2}}
\newcommand{\stala}{\zeta}

\newcommand{\kerov}{\ensuremath{(-1)^l k_1 \cdots k_l
\Sigma_{\seqtwo{k}{l}}}\xspace}
\newcommand{\kerovdiv}{\ensuremath{\Sigma_{\seqtwo{k}{l}}}\xspace}
\newcommand{\kerovn}{\ensuremath{(-1)^l k_1 \cdots k_l
\Sigma_{\seqtwo{k}{r};2n}}\xspace}
\newcommand{\kerovdivnsub}[1]{\ensuremath{\Sigma_{\seqtwo{k}{l};#1}}\xspace}
\newcommand{\mhat}{\ensuremath{\hat{m}}\xspace}
\newcommand{\ptq}{\ensuremath{p \times q}\xspace}
\newcommand{\laginv}{ {\ensuremath{\langle -1 \rangle}}}
\newcommand{\prodofvars}[2]{#1_1  #1_2 \cdots #1_{#2}}
\newcommand{\seq}[2]{\ensuremath{#1_1, #1_2, \ldots, #1_{#2}}}
\newcommand{\coeff}[1]{\ensuremath{[#1]\;}\xspace}
\newcommand{\seqtwo}[2]{#1_1, \ldots, #1_{#2}}
\newcommand{\f}{\ensuremath{\frac}}
\newcommand{\upth}{\mathrm{th}}
\newcommand{\upid}{\mathrm{id}}

\author{Amarpreet Rattan}
\address{Department of Mathematics, Massachusetts Institute of Technology,
Cambridge, MA, 02138, USA}
\email{arattan@math.mit.edu}

\author{Piotr \'Sniady}
\address{Institute of Mathematics,
University of Wroclaw,  \mbox{pl.\ Grunwaldzki~2/4,} 50-384
Wroclaw, Poland} \email{Piotr.Sniady@math.uni.wroc.pl}

\title[Upper bound on characters and a Frobenius formula]{Upper bound on the
characters of the symmetric groups for balanced Young diagrams \\ and a generalized
Frobenius formula}

\DeclareMathOperator{\Tr}{Tr}

\begin{document}

\begin{abstract}
We study asymptotics of an irreducible representation of the symmetric group
$S_n$ corresponding to a balanced Young diagram $\lambda$ (a Young diagram with
at most $C \sqrt{n}$ rows and columns for some fixed constant $C$) 
in the limit as $n$ tends to infinity. We show that there exists a constant
$D$ (which depends only on $C$) with a property that 
$$ |\chi^{\lambda}(\pi)| = 
\left| \frac{\Tr \rho^{\lambda}(\pi)}{\Tr \rho^{\lambda}(e)} \right| \leq 
\left( \frac{D \max(1,\frac{|\pi|^2}{n}  )}{\sqrt{n}} \right)^{|\pi|},$$
where $|\pi|$ denotes the length of a permutation (the minimal number of factors
necessary to write $\pi$ as a product of transpositions).
Our main tool is an analogue of Frobenius character formula which holds
true not only for cycles but for arbitrary permutations. 
\end{abstract}

\maketitle

\section{Introduction}

\subsection{Asymptotics of characters of symmetric groups}

The character of an irreducible representation corresponding to a Young 
diagram $\lambda$ can theoretically be calculated with the help of the Murnaghan--Nakayama
rule. Unfortunately, this combinatorial algorithm quickly becomes intractable when the number of boxes of $\lambda$ tends to infinity. For this reason one is in need of more analytic (possibly approximate) methods which would give meaningful answers asymptotically. 

%

In this article we will concentrate on a scaling in which the number
of rows and columns of a Young diagram $\lambda$ having $n$ 
boxes is bounded from above by
$C\sqrt{n}$ for some fixed
constant $C$.
Such Young diagrams will be called
\emph{balanced}. The investigation of this limit was initiated by Kerov
\cite{Kerov1993transition,Kerov1999differential}
and was further developed by Biane
\cite{Biane1998,Biane2001approximate} who
proved that 
\begin{equation} 
\label{eq:szacowanieBiane}
|\chi^{\lambda}(\pi)| = O(n^{- \frac{|\pi|}{2}}),
\end{equation}
holds true asymptotically for any fixed permutation $\pi$ in the limit $n\to\infty$ for
balanced Young diagrams $\lambda$.  Here, by $|\pi|$ we denote the minimal number of
factors necessary to write $\pi$ as a product of transpositions. The quantity
$\chi^{\lambda}$ appearing in \eqref{eq:szacowanieBiane} is
defined as
$$ \chi^{\lambda}(\pi) = \frac{\Tr \rho^{\lambda}(\pi)}{\Tr
\rho^{\lambda}(e)},$$
where $\rho^\lambda(\pi)$ is the representation of the symmetric group indexed
by $\lambda$ and evaluated at the conjugacy class indexed by $\pi$ and, as
usual, $\rho^\lambda(e)$ is the degree of the character $\rho^\lambda$.
Since this normalized character will appear repeatedly in this article, for
simplicity we will refer to it as the \emph{character} and the name  
\emph{normalized character} will be reserved for the quantity which will be
introduced in \eqref{eq:normalized-characters}.
It should be stressed that the estimate \eqref{eq:szacowanieBiane} was proved by
Biane only asymptotically for $n\gg |\pi|$. 
Biane, in fact, described the asymptotics of \eqref{eq:szacowanieBiane} very
precisely in terms of \emph{free cumulants}, which can be regarded
as certain functionals of the shape of a Young diagram.

The main result of this article is that the inequality
\eqref{eq:szacowanieBiane} holds true for
$|\pi|=O(\sqrt{n})$ and for larger values of $|\pi|$ we still get some meaningful
estimates. 
We state it formally in the following theorem and the proof is given in Section
\ref{sec:asymp}.
\begin{theorem}[Main result]
\label{theo:main}
For every $C>0$ there exists a constant $D$ with the following property. If
$\lambda$ is a Young diagram with $n$ boxes which has at most $C\sqrt{n}$
rows and columns and $\pi\in S_n$ is a permutation then
\begin{equation} 
\label{eq:szacowanie}
|\chi^{\lambda}(\pi)|  <
\left( \frac{D \max(1,\frac{|\pi|^2}{n}  )}{\sqrt{n}} \right)^{|\pi|}.
\end{equation}
\end{theorem}

It should be stressed that for permutations $\pi$
such that $|\pi| \leq o(n^{3/4})$
the inequality \eqref{eq:szacowanie} is significantly stronger
than the best previous estimate for the values of characters on long
permutations, given in \cite{Roichman1996}.  There, the author shows that
$$ |\chi^{\lambda}(\pi)| < q^{|\pi|} $$
holds true for balanced Young diagrams, where $0<q<1$ is a constant.

In the remaining part of this introduction we present the context in which
\eqref{eq:szacowanie} appeared, its applications and the main tools
used in the proof.

\subsection{Generalized Frobenius character formula}
Our strategy in proving estimate \eqref{eq:szacowanie} is to use a generalized
version of Frobenius'
formula for characters, which relates the value of a normalized character evaluated
at some Young diagram with the residue at
infinity of a certain function described by the shape of that Young diagram. The
original Frobenius formula gives the values of the characters only on the cycles and in
Theorem \ref{theo:frobenius} we show its generalization to arbitrary permutations.
This result is interesting in its own and this article seems to be the first in the literature to contain such an expression.

\subsection{Expansion of characters in terms of Boolean cumulants}

For asymptotic problems, it is convenient to encode the shape of a Young diagram
in terms of its \emph{shifted Boolean cumulants}. The generalized Frobenius
formula allows us to express normalized characters as polynomials in shifted Boolean
cumulants by expanding all power series involved. It seems that inequality
\eqref{eq:szacowanie} could be proved directly in this way by crude
estimation of each of the summands, but here we concentrate on a
more elegant approach which we present in the following.

In Theorem \ref{theo:non-negative} we will show that all coefficients in
this expression of characters in terms of shifted Boolean cumulants are
non-negative integers.

\subsection{Bound for characters}

Let Young diagrams $\lambda$, $\mu$ be given with
the shifted Boolean cumulants
$\tilde{B}^{\lambda}_2,\tilde{B}^{\lambda}_3,\dots$ and
$\tilde{B}^{\mu}_2,\tilde{B}^{\mu}_3,\dots$ respectively. Suppose that
\begin{equation}
\label{eq:marzenie}
 |\tilde{B}^{\lambda}_i| \leq \tilde{B}^{\mu}_i 
\end{equation}
holds true for any $i\geq 2$.
The positivity of the coefficients in the expansion of characters in terms of
shifted Boolean cumulants implies that
\begin{equation}
\label{eq:porownanie}
 |\chi^{\lambda}(\pi)| \leq \chi^{\mu}(\pi) 
\end{equation}
holds true for any permutation $\pi$.  Inequality \eqref{eq:porownanie} shows
that in order to prove \eqref{eq:szacowanie} it is enough to prove it for Young
diagrams $\mu$ with sufficiently big shifted Boolean cumulants and we have some
freedom of choosing $\mu$ for which the calculation would take a
particularly nice form. In this article as the reference Young diagram $\mu$ we
take a rectangular Young diagram $p \times q$ since there are very
simple formulae of Stanley for their characters
\cite{StanleyRectangularCharacters}. 

Unfortunately, we encounter some difficulty when trying to follow the
above method, namely for every Young diagram $\mu$ with $n$ boxes the second
shifted Boolean cumulant is given by $\tilde{B}^{\mu}_2=-n<0$ therefore the
inequalities \eqref{eq:marzenie} cannot by fulfilled. A solution to this problem
is to take as $\mu$ a rectangular Young diagram $p\times q$ with $p<0$ and
$q>0$. Clearly, $\mu$ does not make any sense as a Young diagram, nevertheless
we show that for such an object the Stanley character formula still holds true.

We finish the proof of \eqref{eq:szacowanie} with the help of Stanley's character
formulae.

\subsection{Application: asymptotics of Kronecker tensor products and quantum
computations}
\label{subsec:application}
Since this section is not directly connected to the rest of the article we
allow ourselves to be less formal in the following.

The results of this article were motivated by a recent work of Moore and Russell \cite{MooreRussell} who explored quantum algorithms for solving the graph isomorphism problem which are analogous to Kuperberg's algorithm for the dihedral group, namely those which use repeated adaptive tensoring of irreducible representations of the symmetric group $S_n$.
The probability of the success of the algorithm in each iteration depends on the solution to the following problem.
\begin{problem}
\label{problem:decomposition}
Let $\lambda$, $\mu$ be balanced Young diagrams, each having $n$ boxes. Can we
find some upper bound for the (relative) multiplicities in the
decomposition of the Kronecker tensor product $[\lambda]\otimes [\mu]$ into
irreducible components in the limit $n\to\infty$?  In particular, how far is
this tensor product away (in some suitable distance) from the left-regular representation (Plancherel
measure)?
\end{problem}

Moore and Russell conjectured the following stronger version of the estimate \eqref{eq:szacowanie} and proved that
this conjecture would give an upper bound for the multiplicities 
in Problem \ref{problem:decomposition} hence
the expected time of work of the considered quantum algorithm for the graph isomorphism
problem would not be better than the time of work of the best known classical
algorithms.
\begin{conjecture}
\label{conjecture:moore-russell}
For every $C>0$ there exists a constant $D$ with the following property. If
$\lambda$ is a Young diagram with $n$ boxes which has at most $C\sqrt{n}$
rows and columns and $\pi\in S_n$ is a permutation then
\begin{equation} 
\label{eq:szacowanie-lepsze}
|\chi^{\lambda}(\pi)|  <
\left( \frac{D }{\sqrt{n}} \right)^{|\pi|}.
\end{equation}
\end{conjecture}

In a subsequent paper \cite{MooreRussell'Sniady-preprint} we show 
how inequality \eqref{eq:szacowanie} (despite being weaker than  Conjecture
\ref{conjecture:moore-russell}) is sufficient to fill the gap in the work of
Moore and Russell \cite{MooreRussell}.

\subsection{Application: random walks on symmetric groups}
Diaconis and Shahshahani \cite{Diaconis1981} initiated use of the representation
theory in the study of random walks on groups. A typical problem studied in this
context is the speed of convergence of the convolution powers of a given measure
towards the uniform distribution. This method was used with success in numerous
publications (for a review article see
\cite{Diaconis1996}).

Our improved upper bounds for the characters have an immediate application in this context
in the case of a random walk on the symmetric group with steps being random long permutations; the details will be presented in a forthcoming paper \cite{'Sniady2006}.

\subsection{Open problems}

\subsubsection{Optimality of the bound}
The main result of this article, inequality \eqref{eq:szacowanie}, opens many
questions. One problem is to determine how optimal is this bound.  The exact formulae of Stanley \cite{StanleyRectangularCharacters} show that for rectangular Young diagrams and $|\pi|\leq O(\sqrt{n})$ the right-hand side of \eqref{eq:szacowanie} is also a lower bound (with a different value of $D$) and, therefore, this estimate cannot be significantly 
improved.  A possibility of improvement for $|\pi|\geq O(\sqrt{n})$ remains open,
nevertheless some partial results suggest that the estimates as strong as Conjecture
\ref{conjecture:moore-russell} should not be true for general balanced Young diagrams 
and $|\pi|=O(n)$.



\subsubsection{Free cumulants and Kerov polynomials}\label{sec:kerovpolysintro}
In this article we consider the expansion of the normalized characters in
terms of shifted Boolean cumulants while in the asymptotic theory of
representations it is much more common \cite{Biane1998,Biane2001approximate,Sniady04AsymptoticsAndGenus,
Sniady2005GaussuanFluctuationsofYoungdiagrams} to consider analogous expansions
in terms of \emph{free cumulants}. Such expansions of characters in terms of
free cumulants are called \emph{Kerov polynomials}. The main reason for the popularity of free cumulants is that the expansions of characters in terms of free cumulants have a much simpler structure compared to expansions in terms of
Boolean cumulants. For this reason it would be very interesting
to find the analogues of the equalities and estimates presented in this article
in which the shifted Boolean cumulants would be replaced by free cumulants.

For many Young diagrams appearing in the asymptotic theory of representations we
may find much better estimates on free cumulants than on Boolean
cumulants; for
example a typical Young diagram contributing to the Plancherel measure or to a
tensor product of balanced Young diagrams 
has very small free cumulants (except for the second
one) while its Boolean cumulants cannot be bounded better than for a general
balanced Young diagram. It seems plausible that the analogues of
\eqref{eq:szacowanie} obtained by analysis of free cumulants expansions of
characters might give much better estimates for such Young diagrams with
small free cumulants. In particular, if the shapes of the Young diagrams $\mu$
and $\lambda$ are very close to the shape of a typical Young diagram
\cite{LoganShepp,VershikKerov1977} it should be possible to show in this way
that their tensor product $[\mu]\otimes [\lambda]$ in Problem
\ref{problem:decomposition} is indeed very close to the Plancherel measure.


For the method of the proof presented in this article to work with free
cumulants we would have to show that
the coefficients of Kerov polynomials are non-negative integers. This
statement is known as \emph{Kerov's conjecture}. Unfortunately, despite many
results in this direction
\cite{Biane1998,Sniady04AsymptoticsAndGenus,GouldenRattanAccepted,
Biane2005GouldenRattan} Kerov's conjecture remains open until today.

It should also be stressed that usually under the name of Kerov polynomials one
understands the value of $\Sigma_k$ as defined in
\eqref{eq:normalized-characters}, i.e.\ the normalized character evaluated on a
single cycle. For the
purpose of this article such characters are not sufficient and we need to
generalize the definition of Kerov polynomials to cover the expansions for
general permutations.  The problem of positivity of more general Kerov
polynomials was not studied before; in particular it seems possible that in
order to have positive coefficients in free cumulants expansion one should
replace the characters by some kind of cumulants such as the ones considered in
\cite{Sniady2005GaussuanFluctuationsofYoungdiagrams}. We expect that the
generalized Frobenius formula presented in this article will shed some light
into this topic.

In Section \ref{sec:kerovfree} we introduce some of the main ideas concerning
generalized Kerov's polynomials and give some examples.


%

\subsubsection{Combinatorial interpretation}

The positivity of coefficients in the expansion of characters into shifted Boolean cumulants and the (conjectured) positivity of coefficients of Kerov polynomials are rather unexpected and immediately raise questions about some natural combinatorial interpretation to these coefficients. Such interpretations were hinted by Biane \cite{BianeCharacters} but until now there are no concrete results in this direction.

\subsubsection{Approximate factorization of characters}

Biane \cite{Biane2001approximate} proved \emph{approximate factorization
property} for characters of irreducible representations which can be informally
stated as
$$ \chi^{\lambda}(\pi_1 \pi_2) \approx \chi^{\lambda}(\pi_1)
\chi^{\lambda}(\pi_2) $$
for permutations $\pi_1,\pi_2$ with disjoint supports. Using the notation
introduced in \cite{Sniady2005GaussuanFluctuationsofYoungdiagrams} this
property can be formulated formally as the special case $l=2$ of more general
asymptotic estimates on cumulants
\begin{equation}
\label{eq:factorization}
 k_l(\pi_1,\dots,\pi_l) = O(n^{\frac{-|\pi_1|-\cdots-|\pi_l|-2 (l-1)}{2}})
\end{equation}
which hold for balanced Young diagrams $\lambda$ and for
$|\pi_1|,\dots,|\pi_l|\ll n$.
Since the Biane's estimate \eqref{eq:szacowanieBiane} holds true in a more general case
of $|\pi|=O(\sqrt{n})$ (Theorem \ref{theo:main}) it raises a question if estimates
\eqref{eq:factorization} holds true for $|\pi_1|,\dots,|\pi_l|\leq O(\sqrt{n})$.

Biane \cite{Biane1998} and \'Sniady \cite{Sniady2005GaussuanFluctuationsofYoungdiagrams} found not only the asymptotic bound of the right-hand side \eqref{eq:szacowanieBiane} and
\eqref{eq:factorization} for short lengths of the permutations $\pi,\pi_1,\dots,\pi_l$ but also computed explicitly the leading order terms which turn out to be relatively simple functions of free cumulants. It would be very interesting to check if these formulas for
leading order terms hold true without assumptions on
$|\pi|,|\pi_1|,\dots,|\pi_l|$ (our conjecture is that the answer for this
question is negative) and how these leading terms change in the general case.

\subsubsection{Asymptotics of cumulants}
It seems plausible that the generalized Frobenius formula presented in this
article can be used to simplify the calculation of some
of the cumulants in the paper
\cite{Sniady2005GaussuanFluctuationsofYoungdiagrams}.





\section{Continuous Young diagrams and Boolean cumulants}
\label{sec:Young-and-Boolean}


\begin{figure}[tb]
\includegraphics{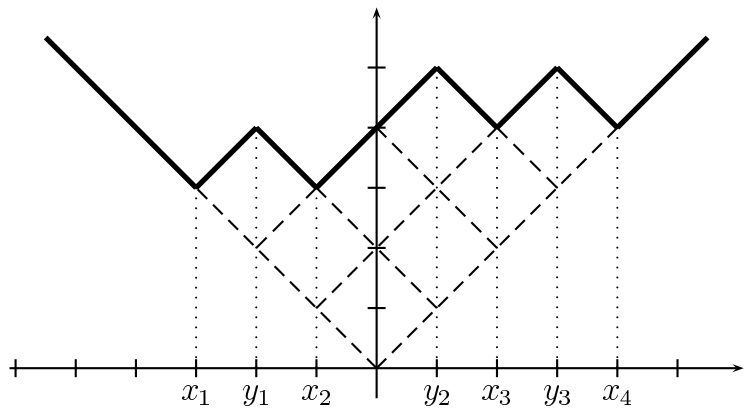}
\caption{Young diagram $(4,3,1)$ drawn according to the Russian convention; 
the ragged bold line represents the profile of a Young diagram. In
this case $(x_1,\dots,x_4)=(-3,-1,2,4)$ and $(y_1,\dots,y_3)=(-2,1,3)$.}
\label{fig:russian}
\end{figure}

Figure \ref{fig:russian} presents our preferred convention for drawing Young
diagrams which is---as opposed to French and English conventions---sometimes
referred to as the \emph{Russian convention}. According to this convention a
Young diagram is represented by its \emph{profile} which is a piecewise affine
function on the real line with the sequence of local minima $x_1<\cdots<x_s$ and
the sequence of local maxima $y_1<\cdots<y_{s-1}$. The modern approach to
representation theory of symmetric groups \cite{OkounkovVershik1996} suggests that
the Russian convention is the most natural one; it originates in the papers of
Kerov
\cite{Kerov1993transition,Kerov1998interlacing}.

To a Young diagram $\lambda$ we associate its \emph{Cauchy transform}  $G(z)$
given by
$$ G(z)= \frac{(z-y_1)\cdots (z-y_{s-1})}{(z-x_1)\cdots (z-x_{s})}. $$
In the following we represent a Young diagram with $n$ boxes as a tuple
$(\lambda_1,\dots,\lambda_n)$ of integers such that $\lambda_1\geq \cdots\geq
\lambda_n\geq 0$, in this case we can also write
$$ G(z)= \frac{(z+n)(z-\lambda_1+1)(z-\lambda_2+2)\cdots (z-\lambda_n+n)}{
(z-\lambda_1)(z-\lambda_2+1)\cdots (z-\lambda_n+n-1)}.$$

In this article we will be concerned mostly with the reciprocal of the Cauchy
transform 
\begin{equation}
 \label{eq:reciprocal}
 H(z)=\frac{1}{G(z)}. 
\end{equation}
For any fixed constant $\stala$ we consider the power series expansion of
$H$ aroung $\zeta$ in decreasing powers of $z$:
$$ H(z+\stala)=\frac{1}{G(z+\stala)}= z +\stala+  \tilde{B}_1 + \tilde{B}_2
z^{-1} + \tilde{B}_3 z^{-2} + \cdots . $$
In fact, one can show that $x_1+\cdots+x_s=y_1+\cdots+y_{s-1}$ therefore
$\tilde{B}_1=0$ and
\begin{equation} 
\label{eq:reciprocal2} 
H(z+\stala)=\frac{1}{G(z+\stala)}= z +\stala+  \tilde{B}_2
z^{-1} + \tilde{B}_3
z^{-2} + \cdots .
\end{equation}
The coefficients $\tilde{B}_2,\tilde{B}_3,\dots$ will be called
\emph{shifted Boolean cumulants}.  This name is coined after \emph{Boolean
cumulants} introduced by Speicher and Woroudi \cite{SpeicherWoroudi} as
coefficients in the expansion of $H$ around $0$:
$$ H(z)=\frac{1}{G(z)}= z -  B_2 z^{-1} -B_3 z^{-2} - \cdots . $$
It follows that if $\stala=0$ shifted Boolean cumulants coincide (up to sign) with the Boolean cumulants of Speicher and Woroudi.  We will call these
cumulants \emph{twisted} Boolean cumulants;  that is, if
\begin{equation}\label{eq:twistbool}
	H(z)=\frac{1}{G(z)}= z +  \hat{B}_2 z^{-1} +\hat{B}_3 z^{-2} - \cdots .
\end{equation}
then $\hat{B}_i$ are called twisted Boolean cumulants, and $\hat{B}_i = -B_i$.  Hereafter, we will always assume that $\zeta\geq 0$.
 
\begin{lemma}
\label{lem:szacowanie-boolowskich}
For any Young diagram whose number of rows and columns are both
smaller than $A$ the
corresponding shifted Boolean cumulant fulfills
$$ |\tilde{B}_k| \leq \frac{[2(A+\zeta)]^k}{2} $$
for each integer $k\geq 2$.
\end{lemma}
\begin{proof}
%
%
Since the sequences of local minima and maxima
interlace $x_1<y_1<x_2<\cdots<y_{s-1}<x_s$, the function $G$ has
positive residues in all of its poles, it is the Cauchy transform of a
certain probability measure $\nu$ called \emph{transition measure} of
the Young diagram $\lambda$:
\begin{equation*}
 G(z)= \int_{\mathbb{R}} \frac{1}{z-x} d\nu(x).
\end{equation*}

The Cauchy transform can be expanded around $\stala$ as a power series in decreasing
powers of $z$:
$$ G(z+\stala)= z^{-1} + M_1 z^{-2}+  M_2 z^{-3} + M_3 z^{-4} + \cdots $$
the coefficients $M_k$ of which are the moments of the shifted transition
measure $$ M_k= \int (x-\stala)^k d\nu(x). $$
Since the support of the transition measure $\nu$ is equal to
$\{x_1,\dots,x_s\}\subset [-A,A]$, we have
\begin{equation}
\label{eq:szacowanie-momenty} 
|M_k| \leq (A+\stala)^k.
\end{equation}
 
We observe that maximal possible absolute value of the shifted Boolean
cumulants
for the sequence $(M_k)$ fulfilling \eqref{eq:szacowanie-momenty} is obtained
for 
\begin{equation}
\label{eq:szacowanie-boolowskich}
 H(z+\zeta)= \frac{1}{z^{-1} - (A+\stala)  z^{-2} - (A+\stala)^2 z^{-3}
-(A+\stala)^3 z^{-4} - \cdots} 
\end{equation}
and a direct calculation of the power series expansion of the right hand side of
\eqref{eq:szacowanie-boolowskich} finishes the proof.
\end{proof}

The constant in the above lemma is not optimal (like all constants in the
following) since in this article we prefer the simplicity of arguments over the
optimality of the constants.

\section{Normalized characters and generalized Frobenius formula}

If $k_1,\dots,k_l$ are positive integers such that $k_1+\dots+k_l=n$ we identify
$(k_1,\dots,k_l)$ with some permutation in $S_n$ with the corresponding cycle
decomposition.  The number of elements in the list $(k_1,\dots,k_l)$ is known as
the \emph{length} of $(k_1,\dots,k_l)$ and is denoted $\ell(k_1,\dots,k_l)$.

Let $\lambda$ be a fixed Young diagram with $n$ boxes. For integers
$k_1,\dots,k_l\geq 1$ we define the \emph{normalized character of a conjugacy
class}
\begin{equation} 
\label{eq:normalized-characters}
\Sigma_{k_1,\dots,k_l} = 
\frac{\Tr
\rho^{\lambda}(k_1,\dots,k_l,1^{n-k_1-\cdots-k_l})}{\Tr \rho^{\lambda}(e)
}
\falling{n}{k_1+\cdots+k_l}, 
\end{equation}
where $\falling{n}{k}=n (n-1) \cdots (n-k+1) $ denotes the falling power. We
use the convention that $\Sigma_{k_1,\dots,k_l} =0$ whenever $k_1+\cdots+k_l>n$.
Normalized characters of conjugacy classes arise naturally
in the study of asymptotics of representations of symmetric groups
\cite{IvanovKerov1999}. 
The following theorem gives a method of computing them.
\begin{theorem}[Generalized Frobenius formula]
\label{theo:frobenius}
For any integers $k_1,\dots,k_l\geq 1$
\begin{multline} 
\label{eq:frobenius}
(-1)^l k_1 \cdots k_l \Sigma_{k_1,\dots,k_l}= \\
[z_1^{-1}] \cdots [z_l^{-1}] \Bigg[ \left(\prod_{1\leq r\leq l} H(z_r) H(z_r-1)
\cdots H(z_r-k_r+1) \right) \\ \prod_{1\leq s<t\leq l}
\frac{(z_s-z_t)(z_s-z_t+k_t-k_s)}{(z_s-z_t-k_s)(z_s-z_t+k_t)} \Bigg]. 
\end{multline}
The right-hand side of \eqref{eq:frobenius} should be understood as follows: we
expand the expression appearing there as a power series in decreasing powers of
$z_l$ with coefficients being functions of $z_1,\dots,z_{l-1}$ and select the
appropriate coefficient. We repeat this procedure with respect to $z_{l-1}$,
$z_{l-2}$,\dots, $z_1$.
\end{theorem}


Before presenting the proof we will introduce some notation.
For sequences of non-negative integers $\mu=(\mu_1,\dots,\mu_n)$ and
$k=(k_1,\dots,k_r)$ let
$\normchar^{\mu}_k$ denote the coefficient of $x^{\mu}=x_1^{\mu_1} \cdots
x_n^{\mu_n}$ in
$$\left(\sum_{1\leq i\leq n} x_i^{k_1}\right) \cdots \left(\sum_{1\leq i\leq
n} x_i^{k_l}\right) \sum_{w\in S_n} (-1)^w x^{w\delta},$$
where $\delta=(n-1,n-2,\dots,1,0)$. 
As seen in \cite{Macdonald}, if $\lambda=(\lambda_1,\dots,\lambda_n)$ is a
Young diagram with $n$ boxes, $\mu=\lambda+\delta$ and $k_1+\cdots+k_l=n$ then
the value of the unnormalized character fulfills
$$ \Tr \rho^{\lambda}(k) = \normchar^{\mu}_k. $$

For $\mu=(\mu_1,\dots,\mu_n)$ we consider the corresponding function
$$ \HH_{\mu}(z)= \frac{(z-\mu_1-1) (z-\mu_2-1)\cdots (z-\mu_n-1) z}{(z-\mu_1)
(z-\mu_2)
\cdots (z-\mu_n)}. $$
In other words, if $\mu=\lambda+\delta$ then 
$$\HH_{\mu}(z+n)=H(z), $$
where $H$ denotes the reciprocal \eqref{eq:reciprocal} of the Cauchy transform
of the Young diagram $\lambda$.

It becomes clear that Theorem \ref{theo:frobenius} will follow from Lemma
\ref{lem:frobenius} below.

\begin{lemma}
\label{lem:frobenius}
For any tuple $\mu=(\mu_1,\dots,\mu_n)$ of non-negative integers and
$k_1,\dots,k_l\geq 1$
\begin{multline*} 
(-1)^l\ \falling{n}{k_1+\cdots+k_l}\ k_1 \cdots k_l\ 
\normchar_{(k_1,\dots,k_l,1^{n-k_1-\cdots-k_l})}^{\mu}
= \\ \normchar_{(1^{n})}^{\mu} \ 
[z_1^{-1}] \cdots [z_l^{-1}] \Bigg[ \left(\prod_{1\leq r\leq l} \HH_{\mu}(z_r)
\HH_{\mu}(z_r-1)
\cdots \HH_{\mu}(z_r-k_r+1) \right) \\ \prod_{1\leq s<t\leq l}
\frac{(z_s-z_t)(z_s-z_t+k_t-k_s)}{(z_s-z_t-k_s)(z_s-z_t+k_t)} \Bigg]. 
\end{multline*}

\end{lemma}
\begin{proof}
We will use induction with respect to $l$. Clearly
$$ \normchar_{(k_1,\dots,k_l,1^{n-k_1-\cdots-k_l})}^{\mu}=
\sum_i
\normchar_{(k_1,\dots,k_{l-1},1^{n-k_1-\cdots-k_{l-1}})}^{\mu^{(i)}}, $$
where $\mu^{(i)}=(\mu_1,\dots,\mu_{i-1},\mu_i-k_l ,\mu_{i+1}, \dots,\mu_n)$
therefore by the induction hypothesis
\begin{multline*}
(-1)^l \falling{n}{k_1+\cdots+k_l} k_1 \cdots k_l
\normchar_{(k_1,\dots,k_l,1^{n-k_1-\cdots-k_l})}^{\mu}=   (-k_l)
\falling{n}{k_l}
\sum_i \normchar_{(1^{n-k_l})}^{\mu^{(i)}} \times \\
[z_1^{-1}] \cdots [z_{l-1}^{-1}] \Bigg[ \left(\prod_{1\leq r\leq l-1}
\HH_{\mu^{(i)}}(z_r)
\HH_{\mu^{(i)}}(z_r-1)
\cdots \HH_{\mu^{(i)}}(z_r-k_r+1) \right) \\ \prod_{1\leq s<t\leq l-1}
\frac{(z_s-z_t)(z_s-z_t+k_t-k_s)}{(z_s-z_t-k_s)(z_s-z_t+k_t)} \Bigg].
\end{multline*}
Since
$$ \HH_{\mu^{(i)}}(z) = \HH_{\mu}(z)
\frac{(z-\mu_i)(z-\mu_i+k_l-1)}{(z-\mu_i-1)(z-\mu_i+k_l)} $$
and (\cite{Macdonald}, page 118)
$$ \falling{n}{k_l} \normchar_{(1^{n-k_l})}^{\mu^{(i)}} =
\normchar_{(1^{n})}^{\mu}  \falling{\mu_i}{r}
\prod_{j\neq i} \frac{\mu_i-\mu_j-k_l}{\mu_i-\mu_j},$$
we have
\begin{multline}
\label{eq:wielkasuma}
(-1)^l \falling{n}{k_1+\cdots+k_l} k_1 \cdots k_l
\normchar_{(k_1,\dots,k_l,1^{n-k_1-\cdots-k_l})}^{\mu}=   (-k_l)
\normchar_{(1^{n})}^{\mu}
  \times \\
[z_1^{-1}] \cdots [z_{l-1}^{-1}] \sum_i \Bigg[ \Bigg( \prod_{j\neq i}
\frac{\mu_i-\mu_j-k_l}{\mu_i-\mu_j} \Bigg) \\
\Bigg(\prod_{1\leq r\leq l-1}
\HH_{\mu}(z_r)
\HH_{\mu}(z_r-1)
\cdots \HH_{\mu}(z_r-k_r+1) \\
\frac{(z_r-\mu_i)(z_r-\mu_i+k_l-k_r)}{(z_r-k_r-\mu_i)(z_r-\mu_i+k_l)} \Bigg) 
\prod_{1\leq s<t\leq l-1}
\frac{(z_s-z_t)(z_s-z_t+k_t-k_s)}{(z_s-z_t-k_s)(z_s-z_t+k_t)} \Bigg].
\end{multline}
It is easy to check that each summand of \eqref{eq:wielkasuma} is equal to the
residue for $z_l=\mu_i$ of the function
\begin{multline}
\label{eq:waznafunkcja}
\normchar_{(1^{n})}^{\mu}
[z_1^{-1}] \cdots [z_{l-1}^{-1}] \Bigg[ \Bigg( \falling{z_l}{k_l}
\prod_{j}
\frac{z_l-\mu_j-k_l}{z_l-\mu_j} \Bigg) \\
\Bigg(\prod_{1\leq r\leq l-1}
\HH_{\mu}(z_r)
\HH_{\mu}(z_r-1)
\cdots \HH_{\mu}(z_r-k_r+1) \\
\frac{(z_r-z_l)(z_r-z_l+k_l-k_r)}{(z_r-k_r-z_l)(z_r-z_l+k_l)} \Bigg) 
\prod_{1\leq s<t\leq l-1}
\frac{(z_s-z_t)(z_s-z_t+k_t-k_s)}{(z_s-z_t-k_s)(z_s-z_t+k_t)} \Bigg]=\\
\normchar_{(1^{n})}^{\mu}
[z_1^{-1}] \cdots [z_{l-1}^{-1}] \Bigg[ 
\Bigg(\prod_{1\leq r\leq l}
\HH_{\mu}(z_r)
\HH_{\mu}(z_r-1)
\cdots \HH_{\mu}(z_r-k_r+1)  \Bigg) \\ 
\prod_{1\leq s<t\leq l}
\frac{(z_s-z_t)(z_s-z_t+k_t-k_s)}{(z_s-z_t-k_s)(z_s-z_t+k_t)} \Bigg].
\end{multline}

Notice that apart from these residues with respect to $z_l$, for each value of
$r\in\{1,\dots,l-1\}$ the series \eqref{eq:waznafunkcja} also has
residues for $z_l=z_r-k_r$ and $z_l=z_r+k_l$.  The residues, however, at these two
points sum to zero,  completing the proof.
\end{proof}

\begin{theorem}
\label{theo:non-negative}
Suppose that positive integers $k_1,\dots,k_l$ are such that 
\begin{equation} 
\label{eq:do-we-have-to-shift}
k_1,\dots,k_l\leq \stala. 
\end{equation}
Then the normalized character 
$(-1)^l \Sigma_{k_1,\dots,k_l}$ is a polynomial in shifted Boolean cumulants
$\tilde{B}_2,\tilde{B}_3,\dots$ with non-negative coefficients. 
\end{theorem}
\begin{proof}
From \eqref{eq:reciprocal2} it follows
that the expansion of 
$$  H(z_r) H(z_r-1) \cdots H(z_r-k_r+1) $$
as a power series in falling powers of $z_r$ and shifted Boolean cumulants has
all non-negative coefficients. 

Without loss of generality we may assume that $k_1\leq\cdots\leq k_l$.
It follows that the factors contributing to \eqref{eq:frobenius} can be written
as
$$ \frac{(z_s-z_t)(z_s-z_t+k_t-k_s)}{(z_s-z_t-k_s)(z_s-z_t+k_t)} =
1+\frac{k_s k_t}{k_s+k_t} \sum_{i\geq 0} \frac{(z_s+k_t)^i-(z_s-k_s)^i
}{z_t^{1+i}}$$
and are, therefore, power series in descending powers of $z_t$, the coefficients
of which are polynomials in $z_s$ with non-negative coefficients.
\end{proof}

It is an interesting question if the condition \eqref{eq:do-we-have-to-shift}
is really necessary; in particular if the above theorem holds true for
$\stala=0$, when the shifted Boolean cumulants coincide 
with the twisted Boolean cumulants introduced in \eqref{eq:twistbool}.  Indeed, the following examples are the
generalized Frobenius formula expanded in terms of twisted Boolean cumulants.
\begin{align*}
	(-1)^{2} \Sigma_{3,2} &=
	{\hat{B}}_{{4}}\hat{B}_{{3}}+13\,\hat{B}_{{2}}\hat{B}_{{3}}+{\hat{B}_{{2}}}^{2}\hat{B}_{{3}}+6\,\hat{B}_{{5}}+18\,\hat{B}_{{3}}\\
	(-1)^{2} \Sigma_{4,2} &=
	\hat{B}_{{5}}\hat{B}_{{3}}+32\,\hat{B}_{{2}}+80\,\hat{B}_{{4}}+8\,\hat{B}_{{6}}+17\,{\hat{B}_{{3}}}^{2}+40\,{\hat{B}_{{2}}}^{2}\\	
	&\hspace{4.0cm}+
3\,\hat{B}_{{2}}{\hat{B}_{{3}}}^{2}
+8\,{\hat{B}_{{2}}}^{3}+24\,\hat{B}_{{2}}\hat{B}_{{ 4}}
\end{align*}
We conjecture that these expressions are always positive.
\begin{conjecture}
The normalized character $\Sigma_{k_1, \ldots, k_l}$ can be expressed as a polynomial
in twisted Boolean cumulants $\hat{B}_i$ with coefficients being non-negative integers.
\end{conjecture}
\section{Asymptotics of characters}\label{sec:asymp}

In the case of a rectangular Young diagram $\lambda=p\times q$ (i.e.\ $p$
parts, all equal to $q$) we have
\begin{equation}\label{eq:stanrect} 
	H(z)= \frac{(z-q)(z+p)}{z+p-q},
\end{equation}
which corresponds to the sequence of shifted Boolean cumulants
\begin{equation}
\label{eq:booleanrectangular}
\tilde{B}_i= -pq (q-p-\stala)^{i-2} \qquad \text{for }i\geq 2.
\end{equation}

\begin{theorem}[Stanley \cite{StanleyRectangularCharacters}]
Let $\pi\in S_{k_1+\dots+k_l}$ be a permutation of type $(k_1,\dots,k_l)$ and let
$\lambda=p\times q$ be a rectangular Young diagram. Then the corresponding
normalized character is given by
\begin{equation} 
\label{eq:stanley}
\Sigma_{k_1,\dots,k_l}^{p\times q} = (-1)^{k_1+\cdots+k_l} \sum_{\sigma_1\sigma_2=\pi}
p^{\kappa(\sigma_1)} (-q)^{\kappa(\sigma_2)}, 
\end{equation}
where the sum runs over all permutations $\sigma_1,\sigma_2\in S_{k_1+\dots+k_l}$ such that
$\sigma_1 \sigma_2=\pi$ and $\kappa(\sigma)$ denotes the number of cycles of permutation $\sigma$.
\end{theorem}

The following Lemma was proved by F\'eray (private communication).
\begin{lemma}
\label{lem:feray}
Let $\pi,\sigma_1,\sigma_2\in S_n$ be such that $\pi=\sigma_1 \sigma_2$.
There exist permutations $\sigma_1',\sigma_2'\in S_n$ such that
$\pi=\sigma_1' \sigma_2'$, $|\sigma_1'|+|\sigma_2'|=|\sigma_1|+|\sigma_2|$,
$|\sigma_2'|= |\sigma_2'\sigma_2^{-1} |+|\sigma_2|$   and every
cycle of $\sigma'_1$ is contained in some cycle of $\sigma_2'$.
Furthermore, $\kappa(\sigma_2') \leq \kappa(\pi)$.

\end{lemma}
\begin{proof}
If every cycle of $\sigma_1$ is contained in some cycle of $\sigma_2$ then $\sigma_1'=\sigma_1$ and
$\sigma_2'=\sigma_2$ have the required property. Otherwise, there exist $a,b\in\{1,\dots,n\}$
such that $a$ and $b$ belong to the same cycle of $\sigma_1$ but not the same
cycle of $\sigma_2$. We define $\sigma_1'=\sigma_1 (a,b)$, $\sigma_2'=(a,b) \sigma_2$ and we iterate this procedure
if necessary. Notice that $|\sigma_1'|=|\sigma_1|-1$,
$|\sigma_2'|=|\sigma_2|+1$, so this procedure will finish after a finite number of steps. It remains to prove that
$|\sigma_2'| \geq |\sigma_2'\sigma_2^{-1} |+|\sigma_2|$ (the opposite
inequality follows from the triangle inequality): notice that $|\sigma_2'|-|\sigma_2|$ is equal to $k$
(where $k$ is the number of steps after which the procedure has terminated) and 
$\sigma_2'\sigma_2^{-1} $ is a product of $k$ transpositions, hence
$|\sigma_2'\sigma_2^{-1}|\leq k $.

For any $a\in\{1,\dots,n\}$ the elements $\sigma_2'(a)$ and 
$\sigma'_1(\sigma'_2(a))=\pi(a)$ belong to the same cycle of $\sigma_1'$ and,
hence, to the same cycle of
$\sigma_2'$. It follows that the elements $a$ and $\pi(a)$ belong to the same
cycle of $\sigma_2'$ and, therefore, every cycle of $\pi$ is contained in some
cycle of $\sigma_2'$.  The result now follows.
\end{proof}

\begin{lemma}
\label{lem:nie-za-duzo}
For any integers $n\geq 1$ and $i\geq 0$ and for any $\pi\in S_n$
\begin{equation}
\label{eq:iloscpermutacji}
 \# \{ \sigma\in S_n : | \sigma | = i \} \leq \frac{n^{2i}}{i!}.
\end{equation}
\end{lemma}
\begin{proof}
Since every permutation in $S_n$ appears exactly once in the product
$$ [1+ (12) ] [1+ (13)+(23) ] \cdots [1+(1n)+\cdots+(n-1,n)], $$
we have
$$ \sum_i x^i\ \#\{ \sigma\in S_n : | \sigma | = i \} = (1+x) (1+2x) \cdots (1+(n-
1) x).$$
Each of the coefficients of $x^k$ on the right-hand side is bounded from above
by the
corresponding coefficient of $e^{x} e^{2x} \cdots
e^{(n-1)x}=e^{\frac{n(n-1)x}{2}}$, finishing the proof.

%
\end{proof}

\begin{theorem} For any Young diagram $\lambda$ with at most $A$ rows and
columns and for any integers $k_1,\dots,k_l\geq 1$
$$
 | \Sigma_{k_1,\dots,k_l}| < 
\begin{cases} (16e^2 A)^{K+l} & \text{for } K\leq 8A,\\
(4eK)^K (4A)^l & \text{for } K\geq 8A, 
\end{cases}
 $$
where $K=k_1+\cdots+k_l$.
\end{theorem}
\begin{proof}
Before beginning the proof, notice that the Murnaghan-Nakayama rule shows that we can
assume $k_1,\dots,k_l\leq 2A$; otherwise $\Sigma_{k_1,\dots,k_l}=0$ holds
trivially. By setting $\stala=2A$, the assumptions of Theorem \ref{theo:non-negative}
are satisfied.

On a purely formal level $\Sigma_{k_1,\dots,k_l}$ is a
polynomial in variables $\tilde{B}_2,\tilde{B}_3,\dots$ and, therefore, it is
well-defined for shifted Boolean cumulants given by
\eqref{eq:booleanrectangular}; in this way we define
$\Sigma_{k_1,\dots,k_l}^{p\times q}$ as a polynomial in $p$ and $q$.

From Lemma \ref{lem:szacowanie-boolowskich} it follows that
\begin{equation} 
\label{eq:dasieoszacowaczgory}
|\tilde{B}^{\lambda}_i| \leq \tilde{B}^{p\times q}_i ,
\end{equation}
where $\tilde{B}^{\lambda}_i$ denotes the shifted Boolean cumulant of $\lambda$
and $\tilde{B}^{p\times q}_i$ is the shifted Boolean cumulant given by
\eqref{eq:booleanrectangular} for 
\begin{equation}
\label{eq:specialchoice}
p=-4A, \qquad q=4A.
\end{equation}
Theorem \ref{theo:non-negative} together with \eqref{eq:dasieoszacowaczgory}
imply that
$$ |\Sigma_{k_1,\dots,k_l}^{\lambda}| \leq |\Sigma_{k_1,\dots,k_l}^{p\times
q}|. $$

We regard $\Sigma^{p\times q}_{k_1,\dots,k_l}$ as a polynomial in variables
$p,q$ and we know that \eqref{eq:stanley}
holds true whenever $p,q$ are positive integers; it follows that
\eqref{eq:stanley} holds as equality between polynomials in $p,q$.
In particular, \eqref{eq:stanley} holds true for \eqref{eq:specialchoice} and
$$ | \Sigma^{\lambda}_{k_1,\dots,k_l}| \leq  \sum_{\sigma_1 \sigma_2=\pi}
(4A)^{\kappa(\sigma_1)+\kappa(\sigma_2)}.$$

We consider a map which to a pair $(\sigma_1,\sigma_2)$ associates any pair 
$(\sigma_1',\sigma_2')$ as prescribed by Lemma \ref{lem:feray}. For any 
fixed $\sigma_2'$ the
permutations $\sigma_2$ such that $|\sigma_2'|= |\sigma_2'\sigma_2^{-1} |+|\sigma_2|$
can be identified with non-crossing partitions of the cycles of $\sigma_2'$ 
(see \cite[Section 1.3]{bianeproperties}).  It follows that the number of such permutations
$\sigma_2$ is equal to the product of appropriate Catalan numbers and, hence, this product is bounded from above by $4^{K}$, where 
$K=k_1+\cdots+k_l$. Therefore
\begin{multline}
\label{eq:w-sumie-pan-niewiele-umie}
\sum_{\substack{\sigma_1,\sigma_2\in S_K,\\ \sigma_1 \sigma_2=\pi}} (4A)^{\kappa(\sigma_1)+\kappa(\sigma_2)} \leq
4^K \sum_{\substack{\sigma_1',\sigma_2'\in S_K, \\ \sigma'_1 \sigma'_2=\pi, \\ 
\kappa(\sigma_2')\leq \kappa(\pi)}} (4A)^{\kappa(\sigma'_1)+\kappa(\sigma'_2)}
\leq \\
4^K \sum_{\sigma'_1\in S_K} (4A)^{\kappa(\sigma'_1)+l} \leq
4^K \sum_{0\leq i\leq K-1} \frac{K^{2i}}{i!} (4A)^{K-i+l},
\end{multline}
where the last inequality follows from Lemma \ref{lem:nie-za-duzo}.

If $K\leq 8A$ we use the estimate 
$$4^K \sum_{0\leq i\leq K-1} \frac{K^{2i}}{i!} (4A)^{K-i+l}
<
4^K (4A)^{K+l} \exp\left( \frac{K^2}{4A}  \right)\leq
(16 e^2 A)^{K+l}.$$

If $K\geq 8A$ then in the sum on the right-hand-side of
\eqref{eq:w-sumie-pan-niewiele-umie} the quotient of consecutive summands is
greater than $2$.  We can, therefore, bound
the sum by the sum of an appropriate geometric series and
$$
4^K \sum_{0\leq i\leq K-1} \frac{K^{2i}}{i!} (4A)^{K-i+l}<
4^K \frac{K^{2K}}{K!} (4A)^{l} < (4eK)^K (4A)^l,$$
where the second last inequality follows from $K! \geq
\left(\tfrac{K}{e}\right)^K$.



\end{proof}


\begin{theorem}
\label{theo:almost-main}
Let $\lambda$ be a Young diagram with $n$ boxes and at most
$C\sqrt{n}$ rows and columns. Then
%
\begin{equation} 
\label{eq:almost-final}
|\chi^{\lambda}(\pi)| \leq 
\left(\frac{\max\left[ (16 e^3 C)^3 , (32 e^2 C)^2 \frac{ |\pi|^2}{n} \right] }{\sqrt{n}}
\right)^{|\pi|}. 
\end{equation}
\end{theorem}
\begin{proof}
Let $k_1,\dots,k_l,1^{n-K}$ be the cycle decomposition of $\pi$
with $k_1,\dots,k_l\geq 2$ and $K=k_1+\cdots+k_l$. Notice that $|\pi|=K-l$,
$K+l\leq 3 |\pi|$, and $K\leq 2|\pi|$.

The inequality $ \falling{n}{K} \geq \left( \frac{n}{e} \right)^K $ shows that
for $K\leq 8C\sqrt{n}$
\begin{equation}
\label{eq:small-support}
  |\chi^{\lambda}(\pi)| = \frac{|\Sigma_{k_1,\dots,k_l}|}{\falling{n}{K}} \leq
\frac{(16 e^3 C)^{K+l}}{(\sqrt{n})^{K-l}}\leq \left(\frac{(16 e^3 C)^3}{\sqrt{n}}\right)^{|\pi|} 
\end{equation}
and for $K\geq 8C\sqrt{n}$
\begin{multline}
\label{eq:big-support}
 |\chi^{\lambda}(\pi)| = \frac{|\Sigma_{k_1,\dots,k_l}|}{\falling{n}{K}} \leq
\left( \frac{16 e^2 C K}{\sqrt{n}} \right)^{K} \frac{1}{(\sqrt{n})^{|\pi|}} \leq \\
\left( \frac{(16 e^2 C K)^2 }{n^{3/2 }} \right)^{|\pi|}< 
\left( \frac{(32 e^2 C |\pi|)^2 }{n^{3/2 }} \right)^{|\pi|}.
\end{multline}
\end{proof}

We are now ready to prove our main theorem.

\vspace{0.3cm}
\noindent {\bf Proof of Theorem \ref{theo:main}}.
The inequality \eqref{eq:almost-final} shows that Theorem \ref{theo:main} holds true
for $D=\max\left[(16 e^3 C)^3 ,     (32 e^2 C)^2\right]$.\hfill $\Box$

\vspace{0.3cm}

The following result will be useful in the study of some quantum algorithms \cite{MooreRussell'Sniady-preprint}.

\begin{corollary}
For every $C$ there exist constants $A> 0$ and $B$ such that if 
a Young diagram $\lambda$ with $n$ boxes has at most $C\sqrt{n}$ boxes in each row and column
then
\begin{equation}
\label{eq:moore-russell} 
\sum_{\substack{\pi\in S_n,\\ |\pi|\leq A n^{4/7}}} |\chi^{\lambda}(\pi)|^4 \leq B .
\end{equation}
\end{corollary}
\begin{proof}
We shall use the notation introduced in the proof of Theorem \ref{theo:almost-main}.
By \eqref{eq:small-support} 
the contribution of the permutations $\pi$ with support size $K$ 
less than or equal to $8C\sqrt{n}$ is
bounded from above by
$$ \sum_{i \geq 0} \frac{n^{2i}}{i!} \left(\frac{(16 e^3 C)^3}{\sqrt{n}}\right)^{4i} <
\exp\left[ (16 e^3 C)^{12}\right]. $$ 

By \eqref{eq:big-support}
the contribution of the permutations $\pi$ with support size $K$
greater than or equal to $8C\sqrt{n}$ is
bounded from above by
\begin{multline*} \sum_{0\leq i \leq A n^{4/7}} \frac{n^{2i}}{i!} 
\left( \frac{(32 e^2 C i)^2 }{n^{3/2 }} \right)^{4i}<
\sum_{0\leq i \leq A n^{4/7}} 
\frac{n^{2i} {e^i}}{i^i} 
\left( \frac{(32 e^2 C i)^2 }{n^{3/2 }} \right)^{4i}= \\
\sum_{0\leq i \leq A n^{4/7}} 
\left( \frac{e(32 e^2 C)^8 i^7  }{n^4} \right)^{i}<1
\end{multline*}
for a suitably chosen constant $A>0$.

\end{proof}

%
%
%
%
%
%
%

\section{Generalized Kerov Polynomials}\label{sec:kerovfree}

Theorem \ref{theo:frobenius} expresses the normalized character $\kerovdiv$ in terms of
the boolean cumulants $B_i$.  It is, however, desirable
to express normalized characters in terms of free cumulants, as discussed in Section
\ref{sec:kerovpolysintro}.  We call the expressions of normalized characters in terms of
free cumulants {\it generalized Kerov polynomials}.
We use a change of variables so as to allow ourselves to use the Lagrange
Inversion Theorem (see \cite[Theorem 1.2.4]{GouldenJacksonCombenum}).
Set $\phi(z) = zH(z^{-1})$.  Then $\phi(z) = 1 - \sum_{i\geq 1} B_i z^i$.  If $K(z) =
z^{-1} + \sum_{k\geq 1} R_k z^{k-1}$, where the $R_k$ are {\it free cumulants}
then by definition
\begin{equation}\label{eq:defk}
	K(z) = (1/H(z))^\laginv
\end{equation}
where $\laginv$ denotes compositional inverse.  Set $R(z) = zK(z)$.  One can
show easily that \eqref{eq:defk} implies
\begin{equation*}\label{eq:relbr}
	\left(\frac{z}{R(z)}\right)^\laginv = \frac{z}{\phi(z)}
\end{equation*}
and, therefore, by Lagrange Inversion we have
\begin{equation}\label{eq:cumul}
	R_{k+1} = -\frac{1}{k} [z^{k+1}]\; \phi(z)^k
\end{equation}
and
\begin{equation*}
	B_{k+1} = \frac{1}{k} [z^{k+1}]\; R(z)^k.
\end{equation*}
With this notation, Theorem \ref{theo:frobenius} becomes
\begin{multline*}
		(-1)^l \prodofvars{k}{l} \Sigma_{\seq{k}{l}} =
		\coeff{z_1^{k_1+1}} \cdots \coeff{z_l^{k_l+1}}\\
		\prod_{r=1}^l \left((1 - z_r)
		\cdots (1 - (k-1)z_r)
		\phi(z_r) \cdots \phi\left(\frac{z_r}{1 - (k_r
		-1)z_r}\right)\right)\\
		\prod_{1 \leq s < t \leq l} \frac{(z_t -
		z_s)(z_t - z_s + (k_t - k_s) z_s z_t)}{(z_t - z_s - k_t z_s
		z_t)(z_t - z_s + k_s z_s z_t)}.
\end{multline*}
At this point it may be useful to give some examples.  In \cite[Page 2]{BianeCharacters},
one can find example of the original Kerov polynomials.  The main conjecture
concerning the original Kerov polynomials is that they are positive in free
cumulants.  Some progress has been made concerning this conjecture in
\cite{BianeCharacters, GouldenRattanAccepted, Sniady04AsymptoticsAndGenus}, but the conjecture remains open.
The following are two examples of generalized Kerov polynomials.
\begin{align}
	\Sigma_{3,2} &=R_{{3}}R_{{4}} -5R_{{2}}R_{{3}}-6R_{{5}}-18R_{{3}}\notag\\
	\Sigma_{2,2,2}  &= R_3^3  - 12 R_3 R_4 + 58 R_3 R_2
	\label{eq:exgenkerov}
     + 40 R_5 + 80 R_3  - 6 R_3 R_2^2 
\end{align}
We use a grading on generalized Kerov polynomials as the one used in \cite{GouldenRattanAccepted}.  Namely, we consider the new
series
\begin{align*}
	\Sigma_{\seq{k}{l};2n}  &= [u^{\left(\sum_i k_i +
	l\right) - 2n}]\Sigma_{\seq{k}{l}}(u)\\
	&= [u^{\left(\sum_i k_i + l\right) - 2n}]\Sigma_{\seq{k}{l}}(R_2 u^2,
	R^3 u^3, \ldots).
\end{align*}
That is, if the {\it weight} of a monomial $R_{i_1}^{j_1} \cdots R_{i_t}^{j_t}$
is given by $\sum_s i_s j_s$, then $\Sigma_{\seq{k}{l};2n}$ consists of the terms
of weight $\sum_i k_i + l -2n$ in $\Sigma_{\seq{k}{l}}$.  Define
\begin{equation}\label{eq:defbigphi}
	\Phi(x,u) = \sum_{i \geq 0} \Phi_i(x) u^i = (1 - ux)
	\phi\left(\frac{x}{1-ux}\right).
\end{equation}
\begin{proposition}\label{thm:genformkerov}
	The following equations hold.
	\newcounter{forkthm}\renewcommand\theforkthm{\thetheorem.\alph{forkthm}}
	\begin{list}{\alph{forkthm}.}{\usecounter{forkthm}}
		\item \label{thm:fork1} For $\seq{k}{l} \geq 1$,
		\begin{equation*}
			\begin{split}	
				\kerov =[z_1^{k_1 + 1}] \cdots
				[z_l^{k_l +
				1}]\; \prod_{r=1}^l \prod_{j=0}^{k_r-1}
				\Phi(z_r,j)\\ \prod_{1 \leq s < t \leq l}
				\frac{(z_t - z_s)(z_t - z_s + (k_t - k_s) z_s
				z_t)}{(z_t - z_s - k_s z_s z_t)(z_t - z_s + k_t
				z_s z_t)}
			\end{split}
		\end{equation*}
		\item \label{thm:fork2} For $\seq{k}{l} \geq 1$,
		\begin{equation*}
			\begin{split}
			(-1)^l k_1 \cdots k_l
\Sigma_{\seqtwo{k}{l};2n} = [u^{2n}][z_1^{k_1 + 1}] \cdots [z_l^{k_l +
				1}]\; \prod_{r=1}^l \prod_{j=0}^{k_r -1}
				\Phi(z_r,ju)\\ \prod_{1 \leq s < t \leq l}
				\frac{(z_t - z_s)(z_t - z_s + (k_t - k_s) z_s
				z_t u)}{(z_t - z_s - k_s z_s z_tu)(z_t - z_s + k_t
				z_s z_tu)} 
			\end{split}
		\end{equation*}
	\end{list}
\end{proposition}
\begin{proof}
	The first part of the proposition is a trivial substitution and the proof
	of the second part is similar to the proof of \cite[Proposition
	4.1]{GouldenRattanAccepted}.
\end{proof}
The following proposition also appears in \cite[Theorem
4.9]{Sniady04AsymptoticsAndGenus}, and is obtained with ease below.
\begin{proposition}\label{thm:topterms}
	There is only one term of highest weight in $\kerovdiv$ and it is
	$R_{k_1} \cdots R_{k_l}$, i.e.
	\begin{equation*}
		\kerovdivnsub{0} = R_{k_1} \cdots R_{k_l}
	\end{equation*}
\end{proposition}
\begin{proof}
	Note that the constant term in Proposition \ref{thm:fork2} can be
	obtained by substituting $u=0$ into the equation.  Doing so we obtain
	\begin{equation*}
		\kerovn = [z_1^{k_1 + 1}] \cdots [z_l^{k_l +
		1}]\; \prod_{r=1}^l \prod_{j=0}^{k_r -1}
		\Phi(z_r,0).
	\end{equation*}
	From the definition of $\Phi(x,u)$ given in \eqref{eq:defbigphi}, we see
	that $\Phi(x,0) = \phi(x)$.  Therefore, we have
	\begin{equation*}
		[z_1^{k_1 + 1}] \cdots [z_l^{k_l + 1}] \prod_{r=1}^l
		\phi(z_r)^{k_r},
	\end{equation*}
	The result now follows from \eqref{eq:cumul}.
\end{proof}
We focus on a special case in what follows.

\subsection{Special case: $l = 2$}

Specializing Proposition \ref{thm:fork1} to the case of two variables we obtain
the formula
\begin{equation*}
	\begin{split}
	\Sigma_{r,s} = \Sigma_r \Sigma_s - [x^{r+1}] [y^{s+1}]
	\prod_{j=0}^{r-1} \Phi(x,j) \prod_{j=0}^{s-1} \Phi(y,j)\\
	\cdot \frac{(xy)^2}{(y - x -rxy)(x-y -sxy)}
\end{split}
\end{equation*}
and from Proposition \ref{thm:fork2} it follows
\begin{equation*}
	\begin{split}
		\Sigma_{r,s}(u) = \Sigma_r(u) \Sigma_s(u) - [x^{r+1}] [y^{s+1}]
		\prod_{j=0}^{r-1} \Phi(x,ju) \prod_{j=0}^{s-1} \Phi(y,ju)\\
		\cdot \frac{(xy)^2 u^2}{(y - x -rxyu)(x-y -sxyu)}
	\end{split}
\end{equation*}
From numerical evidence, we have the following conjecture.
\begin{conjecture}\label{thm:conj}
	The following is positive in free cumulants.
\begin{equation*}	
	[x^{r+1}] [y^{s+1}]
	\prod_{j=0}^{r-1} \Phi(x,j) \prod_{j=0}^{s-1} \Phi(y,j)
	\frac{(xy)^2}{(y - x -rxy)(x-y -sxy)}
\end{equation*}
\end{conjecture}
As noted earlier positivity of the original Kerov polynomials is still open and
it appears that the positivity in Conjecture \ref{thm:conj} is, likewise,
difficult to prove (in fact, the expression in Conjecture \ref{thm:conj} is more
complicated).  Below we give an explicit form for the case $l=2$.  We, in
fact, also assume the both parts are the same, that is if our two parts are
$k_1$ and $k_2$, then $k_1 = k_2$.  We are not, however, able to show positivity
from our explicit expression.

Before giving our explicit expression, we prove a lemma.  In what follows, for a
ring $R$, the ring $R[ [\lambda]]_1$ is the ring of power series with invertible
constant term.
\begin{lemma}
	Suppose that $\phi(\lambda)  \in R[ [\lambda]]_1$, $\omega = t
	\phi(\omega)$ and $F(\lambda) = \lambda^{-k} G(\lambda)$ with
	$G(\lambda) \in R[ [\lambda]]_1$.
	Then,
	\begin{equation*}
		\sum_{n \geq -k} a_n t^n = \f{F(\omega)}{1 - t \phi^\prime
		(\omega)},\;\;\;\;\;\;\;\;\; \mathrm{ with }\;\; a_n = [\lambda^n] F(\lambda)
		\phi^n(\lambda) 
	\end{equation*}
	\label{thm:lagmod}
\end{lemma}
\begin{proof}
	Set $L(\lambda) = \f{G(\lambda)}{\phi^k(\lambda)}$.  Since
	$\phi(\lambda), G(\lambda) \in R[ [\lambda]]_1$, the series $L(\lambda)$
	is a formal power series.  Thus, we have by \cite[Theorem 1.2.4, Part
	2]{GouldenJacksonCombenum}
	\begin{equation*}
		\sum_{n \geq 0} c_n t^n = \f{L(\omega)}{1 - t
		\phi^\prime(\omega)},\;\;\;\;\;\;\;\;\; \mathrm{ with }\;\; c_n
		= [\lambda^n] L(\lambda) \phi^n(\lambda)
	\end{equation*}
	However, we have
	\begin{align*}
		\sum_{n \geq 0} c_n t^n &= \f{L(\omega)}{1 - t
		\phi^\prime(\omega)}\\
		&= \f{\omega^k F(\omega)}{\phi^k(\omega)\left( 1 -
		t\phi^\prime(\omega)\right)}\\
		&=  t^k \f{F(\omega)}{1 - t\phi^\prime(\omega)}.
	\end{align*}
	Thus,
	\begin{equation*}
		\sum_{n \geq 0} c_n t^{n-k} = \f{F(\omega)}{1 -
		t\phi^\prime(\omega)}
	\end{equation*}
	with
	\begin{align*}
		c_n &= [\lambda^n] L(\lambda) \phi^n(\lambda)\\
		&= [\lambda^{n-k}] F(\lambda) \phi^{n-k}(\lambda)
	\end{align*}
	for $n \geq 0$.  Setting for $n \geq -k$
	\begin{equation*}
		a_n = c_{n+k} = [\lambda^n] F(\lambda) \phi^n(\lambda)
	\end{equation*}
	gives the result.
\end{proof}

We now give an explicit form for the expression in Conjecture \ref{thm:conj}.
\begin{equation}
	[u^{2n}][x^{r+1}] [y^{r+1}]
	\prod_{j=0}^{r-1} \Phi(x,ju) \prod_{j=0}^{r-1} \Phi(y,ju)
	\frac{(xy)^2u^2}{( (y - x)^2 -(rxyu)^2)}.
\label{eq:reqs}
\end{equation}
Setting 
\begin{equation*}
	a(x,y) = \frac{1}{r^2} \frac{(rxyu)^2}{(y-x)^2}
\end{equation*}
and noting that
\begin{align*}
	\prod_{j=0}^{r-1} \Phi(x,ju) &= \prod_{j=0}^{r-1}\left( \phi(x) +
	\sum_{i \geq 1} \Phi_i(x) (ju)^i\right)\\
	&= \sum_{\lambda} \mhat_\lambda \Phi_\lambda(x) \phi^{r -
	\ell(\lambda)}(x) u^{|\lambda|}
\end{align*}
where $\mhat_\lambda$ is the monomial symmetric function with the substitution
$x_i = i$ for $i \leq r-1$ and $x_i = 0$ for $i \geq r$ and $\Phi_\lambda(x) =
\prod_i \Phi_{\lambda_i}(x)$ for $\lambda = \lambda_1, \lambda_2, \dots$
we expand \eqref{eq:reqs} in power series in $u$ to obtain
\begin{align}
	&[u^{2n}][x^{r+1}] [y^{r+1}]
	\prod_{j=0}^{r-1} \Phi(x,ju) \prod_{j=0}^{r-1} \Phi(y,ju)
	\frac{1}{r^2} \frac{a(x,y)u^2}{(1 - a(x,y) u^2)}\notag\\
	&= [u^{2n}][x^{r+1}][y^{r+1}] \sum_{\lambda_1,
	\lambda_2} \mhat_{\lambda_1} \mhat_{\lambda_2}\Phi_{\lambda_1}(x) 
	\Phi_{\lambda_2}(y) \phi^{r - \ell(\lambda_1)}(x) \phi^{r -\ell(\lambda_2)}(y)
	\notag\\
	&\hspace{5.0cm}\cdot u^{|\lambda_1|} u^{|\lambda_2|}\frac{1}{r^2} \sum_{i \geq 1} \left(
	a(x,y) u^2 \right)^i\notag\\
	&=[x^{r+1}] [y^{r+1}]  \sum_{i \geq 1} \frac{1}{r^2} a(x,y)^i\notag\\
	& \hspace{2cm}\cdot \sum_{\lambda_1, \lambda_2 \atop |\lambda_1| + |\lambda_2| = 2n-2i}
	\mhat_{\lambda_1} \mhat_{\lambda_2}\Phi_{\lambda_1}(x) \Phi_{\lambda_2}(y)
	\frac{\phi^{r+1}(x)}{\phi^{\ell(\lambda_1)+
	1}(x)}\frac{\phi^{r+2}(y)}{\phi^{\ell(\lambda_2)+ 1}(y)}\notag\\
	&= \frac{1}{r^2} [x^{r+1}] [y^{r+1}] \sum_{i \geq 1}
	a(w_1,w_2)^i\notag\\
	&\hspace{2cm}\cdot \sum_{\lambda_1, \lambda_2 \atop |\lambda_1| + |\lambda_2|
	= 2n-2i} \mhat_{\lambda_1} \mhat_{\lambda_2} \frac{\Phi_{\lambda_1}(w_1)
	\Phi_{\lambda_2}(w_2)}{C(x)C(y) \phi^{\ell(\lambda_1)}(w_1)
	\phi^{\ell(\lambda_2)}(w_2)}\label{eq:finexp2r}
\end{align}
where the second last equality follows from two applications of Lemma
\ref{thm:lagmod} and where
\begin{equation*}
	C(x) = \frac{1}{1 - \sum_{i\geq 2} (i-1) R_i x^i}
\end{equation*}
and
\begin{equation*}
	w_1 = \frac{x}{R(x)},\;\;\;\;\;\;\;\;\;\;\ w_2 = \frac{y}{R(y)}.
\end{equation*}
In \eqref{eq:finexp2r} we have an explicit expression for the terms of weight
$k+2 - 2n$ in $\Sigma_{r,r}$.  For $n=1$, that is for the terms of weight $k$ in
$\Sigma_{r,r}$, we have 
\begin{align}
	\frac{1}{r^2} &[x^{r+1}] [y^{r+1}] \frac{a(w_1,
	w_2)}{C(x)C(y)}\notag\\
	&=\frac{1}{r^2}[x^{r+1}] [y^{r+1}] \frac{\left( 1 - \sum_{i \geq
	2} (i -1) R_i x^i  \right) \left( 1 - \sum_{i \geq
	2} (i -1) R_i y^i  \right) x^2 y^2}{\left( x R(y) - y R(x)
	\right)^2}.\label{eq:expexpress}
\end{align}

For the original Kerov polynomials, some considerable amount of work was needed
to obtain positivity for terms of weight $k-1$ (see
\cite{Sniady04AsymptoticsAndGenus} and \cite{GouldenRattanAccepted}).  We see
that Conjecture \ref{thm:conj} seems at least as difficult.

\subsection{Linear Terms}

Note that Conjecture \ref{thm:conj} can be expressed as \linebreak
$(-1)^{l + 1} k^\upid(\Sigma_{i_1}, \ldots, \Sigma_{i_l})$ is
positive in free cumulants, where \linebreak
$k^\upid(\Sigma_{i_1}, \ldots, \Sigma_{i_l})$ is the $l^\upth$
cumulant (see \cite{Sniady2005GaussuanFluctuationsofYoungdiagrams}).  For example, for $\mu = (r, s)$  we have
\begin{equation*}
	(-1)k^\upid(\Sigma_r,\Sigma_s) = \Sigma_r\Sigma_s - \Sigma_{r,s}
\end{equation*}
and for $\mu = (r,s,t)$ we have
\begin{equation*}
	k^\upid(\Sigma_r, \Sigma_s, \Sigma_t) = \Sigma_{r,s,t} - \Sigma_r \Sigma_{s,t} -
	\Sigma_s\Sigma_{r,t} - \Sigma_t \Sigma_{r,s} + 2\, \Sigma_r \Sigma_s
	\Sigma_t
\end{equation*}
Note that all the linear terms are expressed in $k^\upid(\Sigma_{i_1},
\ldots, \Sigma_{i_l})$;  that is, note
that in the expansion
\begin{align*}
	\Sigma_{r,s,t} &= \Sigma_r \Sigma_{s,t} +
	\Sigma_s\Sigma_{r,t} + \Sigma_t \Sigma_{r,s} - 2\, \Sigma_r \Sigma_s
	\Sigma_t\\
	&+ \coeff{x^{r+1}}\coeff{y^{s+1}}\coeff{z^{t+1}} \prod_{j=0}^{r-1}
	\Phi(x,j) \prod_{j=0}^{s-1} \Phi(y,j) \prod_{j=0}^{t-1}
	\Phi(z,j)\\
	&\cdot \f{rs(xy)^2}{(y-x -rxy)(y-x+sxy)}\cdot \f{rt(xz)^2}{(z-x
	-rxz)(z-x+txz)}\\
	&\cdot \f{st(yz)^2}{(z-y -syz)(z-y+tzy)}
\end{align*}
that the linear terms of $\Sigma_{r,s,t}$ do not occur in the expression $\Sigma_r
\Sigma_{s,t} + \Sigma_s\Sigma_{r,t} + \Sigma_t \Sigma_{r,s} - 2\, \Sigma_r
\Sigma_s \Sigma_t$.
We can prove that the linear terms in \linebreak $(-1)^{l+1}k^\upid(\Sigma_{i_1}, \ldots,
\Sigma_{i_l})$  have
positive coefficients.  
\begin{theorem}\label{thm:lincoeff}
	Suppose that $\mu = k_1 k_2 \cdots k_l \vdash k \leq n$ and let $\omega_\mu$ be a fixed
	element in the conjugacy class $\mu$ in $S_k$.  Then the coefficient
	of $R_{b+1}$ in $(-1)^{l + 1} k^\upid(\Sigma_{i_1}, \ldots, \Sigma_{i_l})$ is
	the number of $k$-cycles $c$ such that $c \cdot \omega_\mu$ has $b$ cycles.
\end{theorem}

\begin{proof}
	For the rectangular shape $\ptq$, the series $H(z)$ is given
	\eqref{eq:stanrect}.
	From \cite[Theorem 1 and Proposition 2]{StanleyRectangularCharacters}
	and \cite[Proposition 4.1]{Rattanpos} the free cumulants
	for this shape are given by
	\begin{equation}\label{eq:freecu}
		R_{b+1} = (-1)^b \sum_{ {u,v \in S_b \atop u \cdot v = (b)}
		\atop \kappa(u) + \kappa(v) = b+1} p^{\kappa(u)}
		(-q)^{\kappa(v)}
	\end{equation}
	where $(b) = (1\; 2\; \cdots b)$ and $\kappa(u)$ is the number of cycles
	in $u$.  However, \cite[Theorem 1]{StanleyRectangularCharacters}
	states
	\begin{equation}\label{eq:stanpolys}
		(-1)^k \Sigma^{\ptq}_{\seqtwo{k}{l}} = \sum_{u,v \in S_k
		\atop u \cdot v = \omega_\mu} p^{\kappa(u)} (-q)^{\kappa(v)}.
	\end{equation}
	Expanding $\Sigma_{\seqtwo{k}{l}}$ in terms of the free cumulants
	given in \eqref{eq:freecu} should give us the same equation as \eqref{eq:stanpolys}.

	The coefficient of $p (-q)^b$ in \eqref{eq:stanpolys} is clearly the
	number of $k$-cycles $c$ such that $c \cdot \omega_\mu$ has $b$ cycles.
	The coefficient of $p (-q)^b$  in \eqref{eq:freecu}, however, is clearly
	the number of $b$-cycles $c$ such that $c \cdot (b)$ has $b$ cycles.
	But the latter number is clearly 1.  Since $R_{b+1}$ accounts for
	all occurrences of $p (-q)^b$ in $(-1)^{l + 1}
	k^\upid(\Sigma_{i_1}, \ldots, \Sigma_{i_l})$, comparing coefficients give us the result.
\end{proof}

We remark this theorem implies that that if $b \geq k$ then the coefficient of
$R_{b+1}$ in $(-1)^{l + 1} k^\upid(\Sigma_{i_1}, \ldots, \Sigma_{i_l})$ is zero.
Also note that when $l$ is even (odd), the coefficient of $R_{b+1}$ is negative
(positive) in $\Sigma_{k_1, \ldots, k_l}$.
Finally, we see that this last theorem also implies that the sum  of the
coefficients of all linear terms of $\kerovdiv$ must be the number of long
cycles in $S_{k_1 + \cdots + k_l}$, namely $(k_1 + \cdots + k_l -1)!$.
One can compare these observations to the normalized characters given in
\eqref{eq:exgenkerov}.

\section*{Acknowledgments}

AR thanks Ian Goulden and John Irving for some useful discussions.  AR is supported
by a \emph{Natural Sciences and Engineering Research Council of Canada}
Postdoctoral Fellowship.

P{\'S} thanks Cristopher Moore, Alexander Russell and Valentin F\'eray 
for very stimulating discussions.  
Research for P{\'S} is supported by the MNiSW research grant P03A 013 30, by the EU Research
Training Network ``QP-Applications", contract HPRN-CT-2002-00279 and  by the EC
Marie Curie Host Fellowship for the Transfer of Knowledge ``Harmonic Analysis,
Nonlinear Analysis and Probability", contract MTKD-CT-2004-013389.

\bibliographystyle{alpha}

\bibliography{biblio2}

\end{document}